\documentclass[11pt]{amsart}

\usepackage{amssymb}
\usepackage{amsxtra}
\usepackage[mathscr]{eucal}
\usepackage{graphics}

\numberwithin{equation}{subsection}

\newtheorem{theorem}{Theorem}[subsection]
\newtheorem{proposition}[theorem]{Proposition}
\newtheorem{lemma}[theorem]{Lemma}
\newtheorem{corollary}[theorem]{Corollary}
\newtheorem{conjecture}[theorem]{Conjecture}

\theoremstyle{definition}

\newtheorem{definition}[theorem]{Definition}
\newtheorem{example}[theorem]{Example}
\newtheorem{examples}[theorem]{Examples}

\newtheorem{remark}[theorem]{Remark}

\renewcommand{\AA}{\mathbb{A}}  %%% Don't need Angstroms
\newcommand{\al}{\alpha}
\newcommand{\alm}{\al_M}
\newcommand{\an}{\text{\normalfont{an}}}

\newcommand{\boldb}{\mathbf{b}}

\newcommand{\BG}{\operatorname{\mathscr{BG}}}

\newcommand{\bm}{\mathbf{m}}
\newcommand{\bn}{\mathbf{n}}

\newcommand{\bu}{\mathbf{u}}
\newcommand{\bv}{\mathbf{v}}
\newcommand{\bw}{\mathbf{w}}
\newcommand{\bx}{\mathbf{x}}

\newcommand{\bz}{\mathbf{z}}
\newcommand{\bzero}{\mathbf{0}}
\newcommand{\CC}{\mathbb{C}}
\newcommand{\charact}{\operatorname{char}}

\newcommand{\D}{\Delta}

\newcommand{\fa}{\mathfrak{a}}

\newcommand{\FF}{\mathbb{F}}

\newcommand{\fp}{\mathfrak{p}}

\newcommand{\G}{\Gamma}
\newcommand{\Gal}{\operatorname{Gal}}
\renewcommand{\ge}{\geqslant}
\newcommand{\GGm}{\mathbb{G}_{\text{m}}}

\newcommand{\hatw}{\widehat{w}}

\newcommand{\kk}{\Bbbk}
\newcommand{\kx}[1]{\kk[x_1^{\pm1},\dots,x_{#1}^{\pm1}]}

\newcommand{\kka}{\goodbar{\kk}}
\newcommand{\kkas}{\kka^{\times}}

\newcommand{\KK}{\mathbb{K}}

\renewcommand{\le}{\leqslant}
\newcommand{\LL}{\mathbb{L}}
\newcommand{\Log}{\operatorname{Log}}
\newcommand{\Max}{\operatorname{Max}}
\newcommand{\N}{\mathsf{N}}
\newcommand{\Nt}{\widetilde{\N}}

\newcommand{\ord}{\operatorname{ord}}

\newcommand{\QQ}{\mathbb{Q}}

\newcommand{\Qpbar}{\overline{\QQ}_p}
\newcommand{\rd}{\RR^d}
\newcommand{\RR}{\mathbb R}

\newcommand{\Rx}{R[\xpm]}
\newcommand{\SA}{\mathscr{A}}
\newcommand{\SB}{\mathscr{B}}
\newcommand{\SC}{\mathscr{C}}
\newcommand{\SF}{\mathscr{F}}
\newcommand{\SN}{\mathscr{N}}
\newcommand{\SNt}{\widetilde{\SN}}
\newcommand{\SO}{\mathscr{O}}
\newcommand{\Spec}{\operatorname{Spec}}
\newcommand{\Sphd}{\mathsf{S}_{d-1}}
\newcommand{\Sd}{\mathsf{S}_{d-1}}

\renewcommand{\SS}{\mathscr{S}}
\newcommand{\ST}{\mathscr{T}}
\newcommand{\SV}{\mathscr{V}}
\newcommand{\SW}{\mathscr{W}}
\newcommand{\td}{T_d}
\newcommand{\ux}{\underline{X}}
\newcommand{\uz}{\underline{Z}}

\newcommand{\val}{\operatorname{val}}
\newcommand{\Var}{\mathsf{V}}

\newcommand{\xpm}{\mathbf{x}^{\pm}}
\newcommand{\zd}{\mathbb{Z}^d}
\newcommand{\Zx}{\ZZ[\xpm]}
\newcommand{\ZZ}{\mathbb{Z}}

\newcommand{\bbl}{[\![} % Laurent polynomials double brackets
\newcommand{\bbr}{]\!]}
\newcommand{\ppl}{(\!(} % Laurent parentheses
\newcommand{\ppr}{)\!)}
\newcommand{\<}{\langle}
\renewcommand{\>}{\rangle}  %%% We don't use \> for tabbing
\newcommand{\aal}{\<\!\<}
\newcommand{\aar}{\>\!\>}

\newsymbol\emptyphi 203F
\renewcommand{\emptyset}{\emptyphi}

\newcommand{\goodbar}{\overline}

\begin{document}

\title[Non-archimedean Amoebas and Tropical Varieties]
  {Non-archimedean Amoebas\\ and Tropical Varieties}

\author[Einsiedler, Kapranov, and Lind]
  {Manfred Einsiedler, Mikhail Kapranov, and
   Douglas Lind}

\address{Manfred Einsiedler, Department of Mathematics,
   Box 354350, University of Washington, Seattle, WA 98195--4350}
\email{einsiedl@math.washington.edu}

\address{Mikhail Kapranov, Department of Mathematics, Yale
  University, New Haven, CT 06520}
\email{mikhail.kapranov@yale.edu}

\address{Douglas Lind, Department of Mathematics,
   Box 354350, University of Washington, Seattle, WA 98195--4350}
\email{lind@math.washington.edu}

\date{\today}

\keywords{amoeba, tropical variety, non-archimedean field,
Bieri-Groves set, total concavity, subdynamics}

\subjclass[2000]{Primary: 14J29, 13A18, 13C15, 13E15; Secondary:
37B05, 22D40}

% % 14J29 Surfaces of general type
% % 13A18 Valuations and their generalizations
% % 13E15 Rings and modules of finite generation or presentation;
% % number of generators
% % 13C15 Dimension theory, depth, related rings (catenary, etc.)
% % 37B05 Transformations and group actions with
% %       special properties (minimality, distality, proximality, etc.)
% % 22D40: Ergodic theory on groups

\thanks{The authors thank Bernd Sturmfels for several crucial
conversations and references, and the American Institute of
Mathematics in Palo Alto for inviting us to the conference
``Amoebas and Tropical Geometry,'' where some of the work was
done. The second author was supported in part by a grant from
NSERC}

\begin{abstract}
   We study the non-archimedean counterpart to the complex amoeba
   of an algebraic variety, and show that it coincides with a
   polyhedral set defined by Bieri and Groves using valuations.
   For hypersurfaces this set is also the tropical variety of the
   defining polynomial.  Using non-archimedean analysis and a
   recent result of Conrad we prove that the amoeba of an
   irreducible variety is connected.  We introduce the notion of
   an adelic amoeba for varieties over global fields, and
   establish a form of the local-global principle for them.  This
   principle is used to explain the calculation of the
   nonexpansive set for a related dynamical system.
\end{abstract}

\maketitle

\section{Amoebas}\label{sec:amoebas}

\subsection{Generalities}\label{subsec:generalities}

Let $\kk$ be a field. Recall \cite[VI.6.1]{Bourbaki} that a norm
(or absolute value) on $\kk$ is a function $a\mapsto |a|$ from
$\kk$ to $\RR_{\ge0}$ such that
\begin{align}
  |a|&=0 \text{ if and only if } a=0,\\
  |ab|&=|a|\,|b|,\\
  |a+b|&\le|a|+|b|. \label{subitem:triangle}
\end{align}
In this paper, unless otherwise specified, we will assume that
$\kk$ is equipped with a nontrivial norm $|\cdot|$ and is complete
with respect to it. If $\LL/\kk$ is an extension of degree $n$,
then the norm on $\kk$ extends to $\LL$ by the formula
$|a|=|N_{\LL/\kk}(a)|^{1/n}$. By $\kka$ we denote a fixed
algebraic closure of $\kk$. Thus a norm on $\kk$ extends to $\kka$
and we have the map
\begin{equation}\label{eqn:Log}
   \Log\colon \bigl( \kka^{\times} \bigr)^d \to\rd, \quad
   (a_1,\dots,a_d)\mapsto (\log
   |a_1|,\dots,\log|a_d|).
\end{equation}
Note that the image of $\Log$ in \eqref{eqn:Log} is dense in
$\rd$ unless the norm is trivial (i.e., unless $|a|=1$ for any
$a\ne0$).

We denote by $\GGm=\mathbb{A}_{\kk}^1\smallsetminus\{0\}$ the
multiplicative group (punctured affine line) over $\kk$, i.e.,
$\GGm=\Spec \kk[x,x^{-1}]$. The algebraic group $\GGm^d$ will be
referred to as the $d$-dimensional algebraic torus over $\kk$.

Let $X\subset\GGm^d$ be a subscheme, i.e., $X=\Spec (\kx{d}/I)$
where $I$ is an ideal.  Then $X(\kka)$, the set of $\kka$-points
of $X$, is a subset in $(\kkas)^d$.

\begin{definition}\label{def:amoeba}
   The \textit{amoeba} of $X$, denoted $\SA(X)$, is the closure
   of $\Log(X(\kka))$ in $\rd$.
\end{definition}

\begin{example}\label{exam:basic}
   Let $\kk=\CC$ with its usual absolute value. Then $\kka=\kk$
   and $\SA(X)$ coincides with $\Log(X(\CC))$ since
   $X(\CC)\subset (\CC^{\times})^d$ is closed and $\Log$ is a
   continuous proper map. The case when $X$ is a hypersurface in
   $(\CC^{\times})^d$ was first considered in \cite{GKZ}. More
   general complex amoebas were studied in \cite{Mik1,Mik2}.
\end{example}

\subsection{Non-archimedean
norms}\label{subsec:non-archimedean-norms} A norm $|\cdot|$ is
called \textit{non-archimedean} if it satisfies
\begin{equation}\label{eqn:non-archimedean-norm}
   |a+b|\le\max\{|a|,|b|\}. \tag{\ref{subitem:triangle}$'$}
\end{equation}
As is well known, non-archimedean norms on $\kk$ are in bijection
with valuations, i.e., maps $v\colon\kk\to\RR\cup\{\infty\}$
satisfying
\begin{align}
  v(a)&=\infty \text{\ \ if and only if $a=0$},\\
  v(ab)&=v(a)+v(b),\\
  v(a+b) &\ge\min\{v(a),v(b)\}.
\end{align}
Explicitly, $v(a)=-\log|a|$. We will use the same letter to denote
the extension of $v$ to $\kka$.

The connection to valuations is the reason why it will be more
natural for us to work with the map
$\val\colon (\kkas)^d\rightarrow\RR^d$ defined by
\begin{displaymath}
   \val(a_1,\ldots,a_d)=\bigl(v(a_1),\ldots,v(a_d)\bigr)=-\Log
   (a_1,\ldots,a_d).
\end{displaymath}

\begin{definition}\label{def:tropical}
   The \textit{tropical variety} $\ST(X)$ of $X$ is the
   closure of $\val(X(\kka))$ in~ $\rd$.
\end{definition}

Clearly the amoeba and the tropical variety satisfy
$\ST(X)=-\SA(X)$.

\begin{examples}
   (a) The field $\kk=\QQ_p$ of $p$-adic numbers has the $p$-adic
   valuation $v_p$ and the $p$-adic norm $|a|_p=p^{-v_p(a)}$,
   with respect to which it is complete. We sometimes write
   $\QQ_{\infty}$ for $\RR$. Every variety defined over the
   rationals therefore has a $p$-adic amoeba for each
   $p\le\infty$, where $p=\infty$ corresponds to the complex
   amoeba. 
   
   (b) Let $\KK$ be any field. The field $\kk=\KK\ppl t\ppr$ of
   formal Laurent series $g(t)=\sum_{j=m}^\infty g_j t^j$, where
   $g_j\in\KK$ and $m\in\ZZ$, has a discrete valuation $\ord$,
   given by $\ord(g(t))=\min\{j\colon g_j\ne0\}\in\ZZ\subset\RR$.
   If $\KK$ is algebraically closed and $\charact(\KK)=0$, then
   the field of Puiseux series
   \begin{displaymath}
      \bigcup_{n\ge1} \,\KK\ppl t^{1/n} \ppr
   \end{displaymath}
   is algebraically closed \cite{Cohn} and thus coincides with
   $\kka$. It has a $\QQ$-valued valuation $\ord$ defined
   similarly to the above.

   Note that the assumption $\charact(\KK)=0$ is necessary. If
   $\charact(\KK)=p$, then the equation $x^p-x=t^{-1}$ has no
   roots in any $\KK\ppl t^{1/n} \ppr$.
   
   (c) Let $\KK$ be an algebraically closed field of any
   characteristic. A \textit{transfinite Puiseux series over
   $\KK$} is a formal sum $g(t)=\sum_{q\in\QQ}g_{q}t^{q}$, where
   $g_{q}\in\KK$ are such that $\text{Supp}(g)=\{q: g_{q}\ne0\}$
   is well-ordered (i.e. every subset of it has a minimal
   element). Such series form a field $\KK\ppl t^{\QQ}\ppr$
   that is always algebraically closed \cite{Ribenboim}.  For
   example, the equation $x^p-x=t^{-1}$ in (b) has the root
   \begin{displaymath}
      x(t)=t^{-1}+t^{-1/p}+t^{-1/p^2}+\dots .
   \end{displaymath}
   One defines a valuation on $\KK\ppl t^{\QQ}\ppr$ by $v(g)=\min
   \{\text{Supp}(g)\}$.
\end{examples}

\subsection{Conventions}\label{subsec:conventions}
We let $\G=v(\kkas)\subset\RR$ denote the valuation group of
$v$. It is a dense divisible subgroup of $\RR$ by our assumption
that $v$ is nontrivial. By a \textit{convex polyhedron} in
$\RR^d$ we mean a subset $\D$ given by a finite system of
affine-linear inequalities
\begin{displaymath}
   \sum_{j=1}^d b_{ij}u_j\ge c_i, \quad i=1,\dots,r.
\end{displaymath}
We say that $\D$ is \textit{$\G$-rational} if the inequalities
above can be chosen such that $b_{ij}\in\ZZ$ and $c_i\in\G$. In
this case $\D\cap\G^d$ is dense in $\D$. Note that we do not
require $\D$ to have full dimension. By a \textit{($\G$-rational)
polyhedral set} $P$ we mean a finite union of ($\G$-rational)
convex polyhedra. We say that $P$ is of \textit{pure dimension
$r$} if all the maximal polyhedra in $P$ have dimension $r$.

\section{Main results and examples}\label{sec:main-results}

\subsection{Amoebas and tropical varieties of polynomials}
\label{subsec:polynomials}
Let
\begin{displaymath}
   f(\bx)=\sum_{\bn\in\zd} a_{\bn}\bx^{\bn},
   \quad \bn=(n_1,\dots,n_d),\ \ \bx^{\bn}=x_1^{n_1}\dots
   x_d^{n_d}
\end{displaymath}
be a Laurent polynomial with coefficients $a_{\bn}\in\kk$, and
let $X=X_f\subset \GGm^d$ be the hypersurface $\{f=0\}$. We
denote $\ST(X_f)$ by $\ST(f)$. For $\bu\in\rd$ define
\begin{equation}\label{eqn:f-tau}
   f^{\tau}(\bu)=\min_{\bn\in\zd} \{v(a_{\bn})+\,\bu\cdot\bn\}.
\end{equation}
Then $f^{\tau}$ is a convex piecewise-linear function on $\rd$
known as the \textit{tropicalization} of $f$ (see \cite{RST} and
\cite{SS} for background). Note that for almost all $\bn$ we have
$a_{\bn}=0$, so $v(a_{\bn})=+\infty$. Therefore $f^{\tau}$ is the
minimum of finitely many affine-linear functions.

\begin{theorem}\label{thm:A(f)-tropical}
   If $f\ne0$ then $\ST(f)$ is equal to the non-differentiability
   locus of~ $f^{\tau}$. In particular, $\ST(f)$ is either empty
   (when $f$ is a monomial), or is a rational polyhedral set of
   pure dimension $d-1$, or is all of $\rd$ (when $f=0$).
\end{theorem}

This can be reformulated as follows. Denote the convex hull of a
set $E\in\rd$ by $\text{Conv(E)}$. Then
\begin{displaymath}
   \SN(f)=\text{Conv}\{\bn: a_{\bn}\ne0\}\subset\rd
\end{displaymath}
is the Newton polytope of $f$. Let
\begin{displaymath}
   \SNt(f) =\text{Conv}\{(\bn,u)\in\zd\times\RR :
   u\ge v(a_{\bn})\}\subset\RR^{d+1}
\end{displaymath}
be the \textit{extended Newton polyhedron} of $f$. Then $\SNt(f)$
projects to $\SN(f)$ by forgetting the last coordinate.  The
following is then an equivalent formulation of Theorem
\ref{thm:A(f)-tropical}.

\begin{corollary}\label{cor:A(f)-tropical}
   \textnormal{(a)} The connected components of $\rd\smallsetminus
   \ST(f)$ are in bijection with vertices of $\Nt(f)$.

   \textnormal{(b)} If $(\bn,v(a_{\bn}))$ is a vertex of $\SNt(f)$,
   then the corresponding component
   $C_{\bn}\subset\rd\smallsetminus\ST(f)$ consists of
   $\bu\in\rd$ such that
   \begin{displaymath}
      \min_{\bm\in\zd}\{v(a_{\bm})+\bu\cdot\bm\}=v(a_{\bn})+\bu\cdot\bn
   \end{displaymath}
   and the minimum on the left side is achieved for exactly one
   $\bm$, namely $\bm=\bn$.

   \textnormal{(c)} The unbounded connected components of
   $\rd\smallsetminus\ST(f)$ correspond to those vertices
   $(\bn,v(a_{\bn}))$ of $\SNt(f)$ that project to a vertex $\bn$
   of $\SN(f)$.
\end{corollary}

Note that $C_{\bn}$ is a convex polyhedral domain.

\begin{example}\label{exam:newton-polygon}
   The case $d=1$ of Theorem \ref{thm:A(f)-tropical}, or,
   equivalently, Corollary \ref{cor:A(f)-tropical} is well known
   (see \cite[Thm.\ 6.4.7]{Gou} or \cite[Exer.\ 
   VI.4.11]{Bourbaki}). If $f(x)=\sum_{j=r}^s a_j x^j$, then the
   Corollary says that the values $v(z_0)$ for roots $z_0\in\kkas$ of
   $f(x)$ are precisely the negatives of the slopes of the
   non-vertical edges of the Newton polygon $\SNt(f)$.
\end{example}

\begin{proof}[Proof of Theorem \ref{thm:A(f)-tropical}]
   Let $\ST=\ST(f)$ and let $\SS$ be the non-differentiability locus
   of $f^{\tau}$. Clearly $\SS$ is a $\G$-rational polyhedral set.

   \begin{lemma}\label{lem:A-in-B}
      $\ST\subset\SS$.
   \end{lemma}
   \begin{proof}
      Since $\SS$ is closed, it is enough to show that
      $\val(X(\kka))\subset\SS$. Let
      $\bu=(u_1,\dots,u_d)\in\val(X(\kka))$, i.e.  $u_i=v(z_i)$
      where $f(z_1,\dots,z_d)=0$. Note that
      $v(a_{\bn}\bz^{\bn})=v(a_{\bn})+\bu\cdot\bn$.  Recall that
      for non-archimedean absolute values, if
      $a_1,\dots,a_r\in\kkas$ with $a_1+\dots +a_r=0$, then there
      are at least two $a_j$ with maximal $|a_j|$.  Since
      $f(\bz)=\sum a_{\bn}\bz^{\bn}=0$, it follows that there are at
      least two terms in the sum whose valuations are both equal
      to $\min\{v(a_{\bn}\bz^{\bn})\}$. This exactly means that
      $f^{\tau}$ is non-differentiable at $\bu$: two
      affine-linear functionals from the set to be minimized
      achieve the same minimal value at $\bu$.
   \end{proof}

   Since $\SS$ is $\G$-rational polyhedral, $\SS\cap\G^d$ is
   dense in $\SS$. To prove that $\SS\subset\ST$ it is therefore
   enough to prove that $\SS\cap\G^d\subset\ST$. Let
   $\bu\in\SS\cap\G^d$. By changing the variables
   \begin{equation}\label{eqn:change-of-variables}
      z_i\mapsto z_i\cdot a_i, \quad a_i\in\kkas, \quad v(a_i)=u_i,
   \end{equation}
   we reduce to the following.

   \begin{lemma}\label{lem:zero-in-B}
      If $\mathbf{0}\in\SS$, then $\mathbf{0}\in\ST$.
   \end{lemma}
   \begin{proof}
      We will find a root $\bz^0$ of $f$ of the form
      \begin{displaymath}
         z^0_i=(t_0)^{b_i}, \quad t_0\in\kkas, \quad v(t_0)=0
      \end{displaymath}
      for an appropriate choice of
      $\boldb=(b_1,\dots,b_d)\in\zd$. Indeed, let
      \begin{displaymath}
         f_{\boldb}(t)=f(t^{b_1},\dots,t^{b_d})
         =\sum_{\bn}a_{\bn}t^{\boldb\cdot\bn}
         \in\kk[t,t^{-1}].
      \end{displaymath}
      The fact that $\mathbf{0}\in\SS$ means that $\SNt(f)$ has a
      face of positive dimension which is horizontal and whose
      height $u$ is minimal. Let $F$ be the maximal face with
      this property. Assume that $\boldb\in\zd$ is generic in the
      following sense: for each edge $[(\bm,u),(\bn,u)]$ of $F$
      we have $\boldb\cdot(\bm-\bn)\ne0$. Then the extended
      Newton polygon $\SNt(f_{\boldb})\subset\RR^{2}$ has a horizontal
      edge of minimal height $u$. By the classical result in
      Example \ref{exam:newton-polygon}, $f_{\boldb}$ has a root $t_0$
      with $v(t_0)=0$.
   \end{proof}
    This completes the proof of Theorem \ref{thm:A(f)-tropical}.
\end{proof}

We now extend the correspondence between components of
$\rd\smallsetminus\ST(f)$ and Laurent series expansions of $1/f$
to the non-archimedean case. This is similar to the known case
$\kk=\CC$ described in \cite[Ch.\ 6, Cor.\ 1.6]{GKZ}.

Let $(\bn,v(a_{\bn}))$ be a vertex of $\SNt(f)$. We then write
\begin{displaymath}
   f(\bx)=a_{\bn}\bx^{\bn}(1+g(\bx)), \quad g(\bx)=\sum_{\bm\ne\bn}
   \frac{a_{\bm}}{a_{\bn}}\,\bx^{\bm-\bn}
\end{displaymath}
and form the Laurent expansion
\begin{equation}\label{eqn:laurent-expansion}
   R_{\bn}(\bx)=\frac{1}{f(\bx)} =a_{\bn}^{-1}\bx^{-\bn}\sum_{n=0}^\infty
   (-1)^ng(\bx)^n
\end{equation}
by using geometric series.

\begin{proposition}\label{prop:laurent-series}
   \textnormal{(a)} $R_{\bn}(\bx)$  is a well-defined Laurent
   series.

   \textnormal{(b)} The domain of convergence of $R_{\bn}(\bx)$ is
   $\val^{-1}(C_{\bn})$.
\end{proposition}

\begin{proof}
   Assume without loss of generality that $\bzero\in
   C_{\bn}$, for otherwise we can make the same change of
   variables as \eqref{eqn:change-of-variables}.
   Assume that $\bu\in C_{\bn}$. Then
   Corollary \ref{cor:A(f)-tropical}(b) implies that
   $v(b_{\bm})+\bu\cdot\bm>0$ for all coefficients
   $b_{\bm}$ of $g$. For $\bu=\bzero$ this shows together with
   completeness of $\kk$ that the sum (possibly infinite)
   defining the coefficients in $R_{\bn}(\bx)$ (as a Laurent series)
   at each $\bx^{\bm}$, $\bm\in \zd$, converges as claimed in (a).
   
   Suppose $\bz$ satisfies $\val(\bz)=\bu\in C_{\bn}$. Then
   $v(b_{\bm}\bz^{\bm})>0$ for all coefficients $b_\bm$ of
   $g$.  Therefore, $R_{\bn}(\bz)$, considered as a series of
   Laurent polynomials but also as a Laurent series, converges.
   In other words the domain of convergence of $R_{\bn}(\bz)$
   contains $\val^{-1}(C_{\bn})$. Note that the domain of
   convergence of the Laurent series $R_{\bn}(\bx)$ is
   $\val^{-1}(D)$ for some convex $D\subset\RR^d$. Furthermore,
   $D$ cannot contain any element of the boundary of $C_{\bn}$
   since $R_{\bn}(\bz)=1/f(\bz)$ for any $\bz$ where $R_{\bn}(\bz)$
   converges. It follows that $D=C_{\bn}$.
\end{proof}

\subsection{Amoebas and Bieri-Groves sets}\label{subsec:BG-sets}
Let $A$ be a commutative ring with unit $1$. Recall from
\cite[VI.3.1, Def.\ 1]{Bourbaki} that a \textit{(ring) valuation}
$w$ on $A$ is a map $w\colon A\to \RR\cup\{\infty\}$ such that
for all $a,b\in A$ we have that
\begin{align}
  w(ab)&=w(a)+w(b),\\
  w(a+b)&\ge\min\{w(a),w(b)\},\\
  w(0)&=\infty \text{ and } w(1)=0.
\end{align}
There may be nonzero elements $a$ in $A$ for which
$w(a)=\infty$. However, $w^{-1}(\infty)$ is easily seen to be a
prime ideal of $A$. Thus if $A$ is a field, then $w$ is a
valuation in the usual sense of Section
\ref{subsec:non-archimedean-norms}.

Let $\kk$, $|\cdot|$, $v$, and $X\subset\GGm^d$ be as before. Let
$A=\kk[X]$ be the coordinate ring of~ $X$, generated by the
coordinate functions $x_1^{\pm1}$, $\dots$, $x_d^{\pm1}$. Define
$\SW(A)$ to be the set of all ring valuations on $A$ extending
$v$ on $\kk$. Let $G=\Gal(\kka/\kk)$. Then there is an embedding
$X(\kka)/G\to\SW(A)$ given by $\bz\mapsto w_{\bz}$, where
$w_{\bz}(f)=v\bigl(f(\bz)\bigr)$. However, $\SW(A)$ is usually
much bigger than $X(\kka)/G$.

\begin{example}\label{exam:bigger}
   Let $d=1$ and $X=\GGm$, so $A=\kk[x,x^{-1}]$. Assume that
   $\G=v(\kkas)\ne\RR$. Fix $u_0\in\RR\smallsetminus \G$. For
   $f(x)=\sum_{j=r}^s a_jx^j\in A$ define
   \begin{displaymath}
      w(f)=\min_{j\in\ZZ} \{v(a_j)+j\,u_0\},
   \end{displaymath}
   which is a ring valuation on $A$ (see Lemma 1 of
   \cite[VI.10.1]{Bourbaki}). But $w$ does not have the form
   $w_z$ for any $z\in\kkas$ since $u_0\notin\G$. Indeed, an easy
   additional argument shows that even with no assumption on $\G$
   this $w$ cannot have the form $w_z$.
\end{example}

Define the map $\beta\colon\SW(A)\to\rd$ by
\begin{displaymath}
   \beta(w)=\bigl(w(x_1),\dots,w(x_d)\bigr)\in\rd.
\end{displaymath}

\begin{definition}\label{def:BG}
   The \textit{Bieri-Groves set} of $X$ is  defined as
   \begin{displaymath}
      \BG(X)=\beta\bigl(\SW(A)\bigr)\subset\rd.
   \end{displaymath}
\end{definition}

\begin{theorem}[Bieri-Groves \cite{BG}]\label{thm:BR-results}
   Let $X\subset\GGm^d$ be an irreducible variety of dimension
   $r$. Then $\BG(X)$ is a $\G$-rational polyhedral set of pure
   dimension $r$.
\end{theorem}

\begin{remark}\label{rem:finite-union}
   Every variety $X$  is a finite union $X=X_1\cup\dots X_s$ of
   irreducible varieties $X_i$. It is elementary to show that
   $\BG(X)=\BG(X_1)\cup\dots\cup\BG(X_s)$ (see \cite[\S2.2]{BG}),
   so that $\BG(X)$ is a closed set in $\rd$. Hence
   $\ST(X)\subset\BG(X)$ by using the valuations
   $w_{\bz}\in\SW(\kk[X])$ with $\bz\in X(\kka)$ as above.
\end{remark}

Our results relating Bieri-Groves sets to amoebas and tropical
varieties are as follows.

\begin{theorem}\label{thm:main}
   If $X\subset\GGm^d$ is an irreducible variety of dimension
   $r$, then $\ST(X)=\BG(X)$. In particular, $\ST(X)$ and
   $\SA(X)$ are $\G$-rational polyhedral sets of pure dimension
   $r$.
\end{theorem}

If $X=X_1\cup\dots\cup X_s$, then clearly
$\ST(X)=\ST(X_1)\cup\dots\cup\ST(X_s)$. By Remark
\ref{rem:finite-union} we obtain the following.

\begin{corollary}\label{cor:arbitrary-varieties}
   Let $X\subset\GGm^d$ be an arbitrary variety. Then $\ST(X)$ and
   $\SA(X)$ are $\Gamma$-rational polyhedral sets.
\end{corollary}

\begin{theorem}\label{thm:connected}
   If $X$ is irreducible, then $\ST(X)$, and hence $\BG(X)$ and
   $\SA(X)$, are connected.
\end{theorem}

\begin{remark}\label{rem:bg}
   Let $I$ be the ideal in $\kk[x_1^{\pm1},\dots,x_d^{\pm1}]$
   defining $X$. Then trivially $\ST(X)\subset \ST(f)$ for every
   $f\in I$. Speyer and Sturmfels \cite[Thm.\ 2.1]{SS} have shown
   that
   \begin{displaymath}
      \ST(X)=\bigcap_{f\in I} \ST(f).
   \end{displaymath}
   Furthermore, they describe in \cite[Cor.\ 2.3]{SS} that the
   intersection can be taken over just those $f$ in a (finite)
   universal Gr\"obner basis for $I$. Hence a tropical variety is
   always the intersection of a finite number of tropical
   hypersurfaces, each of which has an explicit description as a
   $\Gamma$-rational polyhedral set from Theorem
   \ref{thm:A(f)-tropical}. Their approach can be developed into
   an alternative proof of Theorem \ref{thm:main}.
\end{remark}

\subsection{Adelic amoebas}\label{subset:adelic}
Let $\FF$ be a field of one of the following two types:

(a) A number field, i.e., a finite extension of $\QQ$,

(b) A function field, i.e., $\FF=\KK(C)$ is the field of rational
functions on a smooth projective algebraic curve $C$ of a field
$\KK$ (a detailed account of such fields is contained in \cite{BP}).

Two norms on $\FF$ are said to be equivalent if they define the
same topology. Let $S=S(F)$ be the set of equivalence classes of
norms on $\FF$ (inducing the trivial norm on $\KK$ in the case
(b)). It is well known that one can choose a representative
$|\cdot|_p$ for each $p\in S$ such that:
\begin{equation}\label{eqn:product-formula-1}
   \text{for every $a\in\FF^{\times}$ we have $|a|_p=1$ for almost all
   $p\in S$,}
\end{equation}
\begin{equation}\label{eqn:product-formula-2}
   \prod_{p\in S} |a|_p=1  \quad \text{for every $a\in\FF^{\times}$.}
\end{equation}

We assume that such a choice has been made.

Let $X\subset\GGm^d$ be an algebraic variety defined over
$\FF$. For every $p\in S$ we have the completion $\FF_p$ of $\FF$
with respect to $|\cdot|_p$. Thus we can form the amoeba
\begin{equation}\label{eqn:p-adic-amoeba}
   \SA_p(X)=\text{closure}\bigl[\Log\bigl(X(\overline{\FF}_{p})\bigr)
   \bigr]\subset\rd,
\end{equation}
corresponding to $X$ regarded as a variety over $\FF_p$.

\begin{definition}\label{def:adelic-amoeba}
   The \textit{adelic amoeba }of $X$ is the union
   \begin{equation}\label{eqn:adelic-amoeba}
      \SA_{\AA}(X)=\bigcup_{p\in S} \SA_p(X).
   \end{equation}
\end{definition}

\begin{remark}\label{rem:adelic}
   Let $I$ be the ideal in $\FF[x_1^{\pm1},\dots,x_d^{\pm1}]$
   defining $X$. By the result of Speyer and Sturmfels
   \cite[Cor.\ 2.3]{SS} described in Remark~\ref{rem:bg}, there
   is a finite set $\SF\subset I$ such that for every $p$ we have
   that $\SA_p(X)=\bigcap_{f\in\SF} \SA_p(f)$. Furthermore,
   by Theorem \ref{thm:A(f)-tropical}, for each $f\in\SF$ the
   sets $\SA_p(f)$ are equal for almost every~ $p$. Hence the sets
   $\SA_p(X)$ agree for almost every~ $p$, and so the union in
   \eqref{eqn:adelic-amoeba} can therefore be expressed as a
   finite union. This shows that $\SA_{\AA}(X)$ is a closed
   polyhedral set.
\end{remark}

\begin{remark}\label{rem:BG}
   Bieri and Groves explicitly introduced a ``global'' version of
   their sets by taking the union over all non-archimedean $p$
   (see \cite[Thm.\ B]{BG}). However, for certain results (such as
   the following) the union over all $p$ is needed.
\end{remark}

\begin{theorem}\label{thm:ray}
   Let $X\subset\GGm^d$ be a hypersurface defined over
   $\FF$. Assume that $\mathbf{0}\notin\SA_p(X)$ for at least one
   $p$. Then for any nonzero vector $\bv\in\rd$ the open half-line
   $(0,\infty)\cdot \bv$ meets $\SA_{\AA}(X)$.
\end{theorem}

A special case of the above theorem appears in
\cite[Prop.~5.5]{ELMW} motivated by algebraic dynamical systems
(see Section \ref{sec:connections}), where the proof makes use of
the notion of homoclinic points for such actions. An alternative
proof for the case $d=2$ is worked out in \cite{Scott}.

\begin{example}\label{exam:adelic}
   Let $\FF=\QQ$. Then $S$ consists of all prime numbers and
   $\infty$, with $|x|_p=p^{-\ord_p(x)}$ and $|x|_\infty$ the
   usual absolute value.
   
   Let $d=2$ and $X$ be given by the equation $f(x,y)=3+x+y$.
   Theorem \ref{thm:A(f)-tropical} implies that $\SA_3(X)$ is the
   union of the three rays in Figure~ \ref{fig:3+x+y}(a) meeting
   at the point $(-1,-1)$, and for $p\ne 3$ or $\infty$ the set
   $\SA_p(X)$ is the union of the three rays starting at the
   origin shown in Figure~ \ref{fig:3+x+y}(b). Finally,
   $\SA_\infty(X)$ is shown in Figure~ \ref{fig:3+x+y}(c).
   Observe that every open half-line starting at the origin hits
   at least one of these three amoebas, and in this example
   exactly one amoeba, except in the direction of~ $(1,1)$.
   \begin{figure}[htbp]
      \begin{center}
         \scalebox{1}{\includegraphics{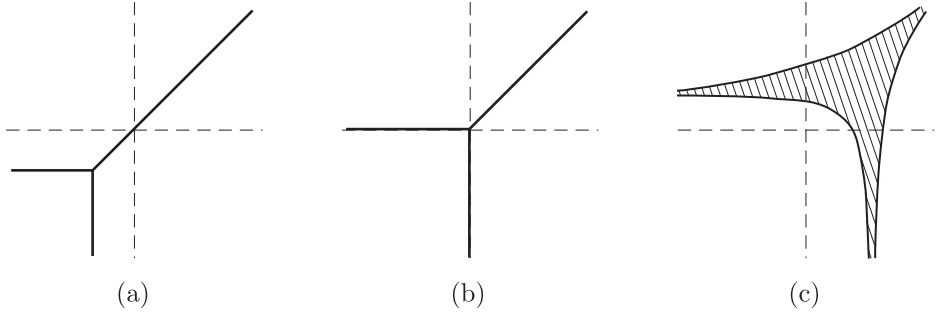}}
         \caption{The amoeba of $3+x+y$ over different fields}
         \label{fig:3+x+y}
      \end{center}
   \end{figure}
\end{example}

We conjecture that a version of this result holds for
lower-dimensional varieties as well.

\begin{conjecture}\label{conj:half-space}
   Let $X\subset\GGm^d$ be irreducible of dimension $r$ such that
   $\SA_{\AA}(X)$ is not contained in any hyperplane. Then for
   any linear subspace $L\subset\rd$ of co-dimension $r$ and any
   relatively open half-space $H$ in $L$ whose boundary contains~
   $\bzero$, we have that $H\cap\SA_{\AA}(X)\ne\emptyset$.
\end{conjecture}

\begin{proof}[Proof of Theorem \ref{thm:ray}]
   Let $f(\bx)=\sum a_{\bn}\bx^{\bn}$ be the equation of $X$, so
   $a_{\bn}\in\FF$. We write $\SA_p(f)$ for $\SA_p(X_f)$ and
   $\ST_p(f)$ for $\ST_p(X_f)$.

   Let $S_{\text{gen}}\subset S$ be the set of all
   non-archimedean norms $p$ for which $|a_{\bn}|_p=1$ for all
   nonzero $a_{\bn}$. Since the number of archimedean norms is
   finite, by \eqref{eqn:product-formula-1} we see that
   $S\smallsetminus S_{\text{gen}}$ is finite, and so
   $S_{\text{gen}}$ is nonempty.

   Let $\SV(f)$ denote the set of vertices of $\SN(f)$. For $p\in
   S_{\text{gen}}$ the extended Newton polytope $\SNt(f)$ is the
   product $\SN(f)\times[0,\infty)$. For each vertex
   $\bn\in\SV(f)$ put
   \begin{displaymath}
      D_{\bn}=\bigl\{\, \bu\in\rd: \bu\cdot\bn>\bu\cdot\bw
      \text{ for all $\bw\in\SN(f)\smallsetminus\{\bn\}$}
      \,\bigr\},
   \end{displaymath}
   which is the \textit{normal cone} to $\SN(f)$ at $\bn$.

   By Corollary \ref{cor:A(f)-tropical}(c), each vertex
   $\bn\in\SV(f)$ corresponds to an unbounded component
   $C_{\bn}^{(p)}$ of $\rd\smallsetminus\ST_p(f)$. Putting
   $D_{\bn}^{(p)}= -C_{\bn}^{(p)}$, using $\SA_p(f)$ instead of
   $\ST_p(f)$ and using $\Log$ instead of $\val$, we see by
   Proposition \ref{prop:laurent-series} that the Laurent series
   $R_{\bn}^{(p)}(\bx)$ in \eqref{eqn:laurent-expansion} for
   $1/f(\bx)$ has domain of convergence
   $\Log^{-1}(D_{\bn}^{(p)})$. Note also that the geometric
   series for $R_{\bn}^{(p)}(\bx)$ produces Laurent coefficients
   that are only finite sums, so that $R_{\bn}^{(p)}(\bx)$ has
   coefficients in $\FF$. Hence $R_{\bn}^{(p)}(\bx)=R_{\bn}(\bx)$
   is independent of $p$. From Corollary
   \ref{cor:A(f)-tropical}(b) we see that $D_{\bn}^{(p)}=D_{\bn}$
   for all $p\in S_{\text{gen}}$. In this case, all components of
   $\rd\smallsetminus \SA_p(f)$ are unbounded, and $\SA_p(f)$ is
   the complement of $\bigcup_{\bn\in\SV(f)} D_{\bn}$ whenever
   $p\in S_{\text{gen}}$.

   Fix $\bv\ne\bzero$, and suppose that $(0,\infty)\cdot\bv$ does
   not meet $\SA_{\AA}(f)$. If $\bv\notin D_{\bn}$ for all
   $\bn\in\SV(f)$, then $\bv\in\SA_p(f)$ for all $p\in
   S_{\text{gen}}$, and we are done. Hence we may assume that
   $\bv\in D_{\bn}$ for some $\bn\in\SV(f)$.

   By hypothesis, there is a $p_0\in S$ so that
   $\bzero\notin\SA_{p_0}(f)$. Since
   \begin{displaymath}
      \bigl[ (0,\infty)\cdot\bv\bigr]\cap\SA_{p_0}(f)
      \subset  \bigl[ (0,\infty)\cdot\bv\bigr]\cap\SA_{\AA}(f)
      =\emptyset,
   \end{displaymath}
   it follows that $[0,\infty)\cdot\bv\subset
   D_{\bn}^{(p_0)}$. Now $D_{\bn}^{(p_0)}$ is convex and open, so
   there is an $\epsilon>0$ so that $R_{\bn}(\bx)=\sum
   b_{\bm}\bx^{\bm}$ converges for all
   $\bx\in\Log^{-1}\bigl[(-\epsilon,\epsilon)\cdot\bv\bigr]$. It
   follows there are $\theta>1$ and $c>0$ so that
   \begin{displaymath}
      |b_{\bm}|_{p_0}<c\,\theta^{-\epsilon|\bm\cdot\bv|}.
   \end{displaymath}
   Also, $|b_{\bm}|_p\le1$ for all $p\in S_{\text{gen}}$. Let
   $r=|S\smallsetminus S_{\text{gen}}|$. By the product formula
   \eqref{eqn:product-formula-2}, there is a $q\in S$ for which
   \begin{displaymath}
      |b_{\bm}|_q\ge c^{1/r}\,\theta^{(\epsilon/r)|\bm\cdot\bv|}
   \end{displaymath}
   for infinitely many $\bm$. This implies that the series
   $\sum b_{\bm}\bx^{\bm}$ does not converge with respect to
   $|\cdot|_q$ for any point in
   $\Log^{-1}\bigl[(-\epsilon/r,\epsilon/r)\cdot\bv)\bigr]$.
   Hence the unbounded component $D_{\bn}^{(q)}$ does not meet
   $(-\epsilon/r,\epsilon/r)\cdot\bv$, and so $(0,\infty)\cdot\bv$
   must meet $\SA_q(f)$.
\end{proof}

\begin{remark}\label{rem:transcendental}
   The proof of Theorem \ref{thm:ray} depends crucially on the
   observation that the Laurent expansions for unbounded
   components have coefficients in the ground field. This can
   fail for bounded components. For example, let $\FF=\QQ$ and 
   $f(x,y)=4-x-y-x^{-1}y^{-1}$. It is shown in \cite[Exam.\
   5.8]{LS} that the constant term in the Laurent expansion of
   $1/f$ for the bounded component containing the origin is 
   \begin{displaymath}
      \frac{1}{4}\sum_{n=0}^\infty \frac{(3n)!}{(n!)^3}\,
      4^{-3n},
   \end{displaymath}
   and is the value of a hypergeometric function known to be
   transcendental.
\end{remark}

\section{Non-archimedean analysis and proofs of main
   results}\label{sec:proofs}

\subsection{Affinoid algebras and affinoid
     varieties}\label{subsec:affinoids}
Let 
\begin{displaymath}
   \td=\kk\aal x_1,\dots,x_d\aar\subset\kk\bbl
   x_1,\dots,x_d\bbr
\end{displaymath}
be the set of formal series
$f(\bx)=\sum_{\bn\in\ZZ_+^d} a_{\bn}\bx^{\bn}$ such that
$|a_{\bn}|\to0$ as $\|\bn\|\to\infty$, where
$\|\bn\|=|n_1|+\dots+|n_d|$. 
This set is a $\kk$-algebra
called the \textit{Tate algebra}. Throughout this section $G$
denotes the Galois group $\Gal(\kka/\kk)$.

\begin{proposition}\label{prop:tate-algebra}
   \textnormal{(a)} The ring $\td$ is Noetherian.

   (b) The maximal ideals of $\td$ are in bijection with the points
   of the unit polydisk
   \begin{displaymath}
      \{\bz\in\kka^d:|z_j|\le1 \text{ for $1\le j\le d$}\}
   \end{displaymath}
   modulo the action of the Galois group $G$.
   Explicitly, if $\bz=(z_1,\dots,z_d)$ is
   such a point whose coordinates $z_j$ generate the
   finite extension $\LL$ of $\kk$,
   then for any $f\in\td$ the series $f(\bz)$ converges in $\LL$,
   yielding a surjective homomorphism $\td\to\LL$.
\end{proposition}

\begin{proof}
   For (a) see Theorem 1 of Section 5.2.6 of \cite{BGR}, and for
   (b) see Proposition 1 of Section 7.1.1. of the same reference.
\end{proof}

For $f(\bx)=\sum a_{\bn}\bx^{\bn}\in\td$ we set
$|f|=\max\{|a_{\bn}|\}$. This makes $\td$ into a $\kk$-Banach
algebra, i.e., a complete non-archimedean normed algebra with
norm extending that on $\kk$.

\begin{definition}\label{def:affinoid-algebra}
   An \textit{affinoid algebra} is a $\kk$-Banach algebra $A$
   admitting a continuous epimorphism $\td\to A$ for some $d$.
\end{definition}

\begin{examples}\label{exam:affinoid}
   (a) All ideals in $\td$ are closed (see Corollary 2 of
   \cite[\S5.2.7]{BGR}). If $I\subset\td$ is an ideal, then
   $A=\td/I$ has the residue norm $|f+I|=\inf\{|f+g|:g\in I\}$,
   which makes it into an affinoid algebra via the projection
   $\td\to A$. Up to replacing a norm with an equivalent one
   (giving the same topology) this is a general form of an
   affinoid algebra.

   (b) If $B$ is an affinoid algebra, $A$ a $\kk$-Banach algebra,
   and $\varphi\colon B\to A$ a continuous homomorphism such that
   $A$ is finitely generated as a $B$-module, then $A$ is affinoid.
\end{examples}

\begin{example}\label{exam:laurent-affinoid}
   Let $\D\subset\rd$ be a bounded $\G$-rational
   polyhedron. We define $\kk\aal\val^{-1}(\D)\aar$ to consist of
   formal Laurent series $f(\bx)=\sum_{\bn\in\zd} a_{\bn}\bx^{\bn}$
   satisfying the condition
   \begin{displaymath}
      \lim_{\|\bn\|\to\infty} \{ v(a_{\bn})+\bn\cdot\bu\}=\infty
      \text{ for
      all $\bu\in\D$}.
   \end{displaymath}
   One sees directly that $\kk\aal\val^{-1}(\D)\aar$ is a ring
   and the norm corresponding to the valuation
   \begin{displaymath}
      v(f)=\inf\{\,v(a_{\bn})+\bn\cdot\bu:\bn\in\zd,\bu\in\D\,\}
   \end{displaymath}
   makes it into a $\kk$-Banach algebra.
\end{example}

Note that for any $\bz\in(\kka^\times)^d$ such that $\val(\bz)\in\D$
the series $f(\bz)$ converges for any
$f\in\kk\aal\val^{-1}(\D)\aar$.

\begin{proposition}\label{prop:affinoid-proof}
   $\kk\aal\val^{-1}(\D)\aar$ is an affinoid algebra.
\end{proposition}

We will call such algebras \textit{polyhedral affinoid algebras}.

\begin{proof}
   Write the inequalities defining $\D$ in the form
      \begin{displaymath}
         \sum_{j=1}^d b_{ij} u_j\ge c_i, \quad i=1,\dots,r
      \end{displaymath}
      with $b_{ij}\in\ZZ$ and $c_i\in\G=v(\kka^\times)$. Without
      loss of generality we can assume $c_i\in v(\kk^\times)$ by
      taking, if necessary, integer multiples of the
      inequalities. Let $\boldb_i=(b_{i1},\dots,b_{id})\in\zd$,
      and choose $w_i\in\kk^\times$ with $c_i=v(w_i)$. Then the
      condition $\val(\bx)\in\D$ can be rewritten as
      $|\bx^{\boldb_i}/w_i|\le1$ for $i=1,\dots,r$. Thus there is
      a continuous homomorphism
      \begin{displaymath}
         T_r=\kk\aal y_1,\dots,y_r\aar
         \stackrel{\varphi}{\longrightarrow}
         \kk\aal\val^{-1}(\D)\aar, \quad y_i\mapsto \bx^{\boldb_i}/w_i.
      \end{displaymath}
      By Example \ref{exam:affinoid}(b) it is enough to prove
      that $\varphi$ is finite. Let $S\subset\zd$ be the
      semigroup with $\mathbf{0}$ generated by $\boldb_i$,
      $i=1,\dots,r$, and $C=\text{Conv}(S)$ its convex
      hull. Denote by $\kk\bbl S\bbr$ the set of all formal
      Laurent series $\sum_{\bn\in S}a_{\bn}\bx^{\bn}$, and by
      $\kk\aal\val^{-1}(\D)\aar_S$ the intersection
      $\kk\aal\val^{-1}(\D)\aar\cap\kk\bbl S\bbr$. Our statement
      reduces to the following statements.
      \begin{lemma}\label{lem:affinoid}
         \textnormal{(a)}
         $\text{\textnormal{Im}}(\varphi)=\kk\aal\val^{-1}(\D)\aar_S$.

         \textnormal{(b)} $\kk\aal\val^{-1}(\D)\aar$ is finite over
         $\kk\aal\val^{-1}(\D)\aar_S$.
      \end{lemma}
      \begin{proof}
         (a) Let $f(\bx)=\sum_{\bn\in S}a_{\bn}\bx^{\bn}$ lie in
         $\kk\aal\val^{-1}(\D)\aar_S$. For any $\bn\in S$ choose
         $\bm(\bn)=(m_1(\bn),\dots,m_r(\bn))\in\ZZ_+^r$ such that
         $m_1(\bn)\boldb_1+\dots+m_r(\bn)\boldb_r=\bn$. Then the series
         $g(\mathbf{y})=\sum_{\bn\in S}a_{\bn}
         \mathbf{y}^{\bm(\bn)}$ lies in $T_r$
         and $\varphi(g)=f$.

         (b) Let $\kk[S]$ and $k[C\cap\zd]$ be the semigroup
         algebras of $S$ and of $C\cap\zd$. Then it is well known
         that $\kk[C\cap\zd]$ is finite over $\kk[S]$. A system
         of module generators of $\kk[C\cap\zd]$ over $\kk[S]$
         will be a system of module generators of
         $\kk\aal\val^{-1}(\D)\aar$ over $\kk\aal\val^{-1}(\D)\aar_S$.
      \end{proof}
   This completes the proof of Proposition \ref{prop:affinoid-proof}
\end{proof}

\begin{example}\label{exam:delta=0}
   A particular case of a polyhedral affinoid algebra is
   $\D=\{\mathbf{0}\}$. The algebra
   $\kk\aal\val^{-1}(\mathbf{0})\aar$ consists of Laurent series
   $\sum_{\bn\in\zd} a_{\bn}\bx^{\bn}$ with $|a_{\bn}|\to0$ as
   $\|\bn\|\to\infty$. It is
   the quotient of $T_{2d}=\kk\aal z_1,w_1,\dots,z_d,w_d\aar$ by
   the ideal generated by the elements $z_iw_i-1$, $i=1,\dots,d$.
\end{example}

For an affinoid algebra $A$ we denote by $\Max(A)$ the set of its
maximal ideals. Recall that $G=\Gal(\kka/\kk)$.

\begin{proposition}\label{prop:maximal}
   \textnormal{(a)} $\Max(A)\ne\emptyset$ unless $A=0$.

   \textnormal{(b)} If $A=\kk\aal z_1,\dots,z_d\aar/
   \<f_1,\dots,f_r\>$, then
   $\Max(A)$ is identified with the set of
   $G$-orbits on
   \begin{displaymath}
      \{\bz\in\kka^d:
      |z_j|\le1 \text{ for $1\le j\le d$ and $f_i(\bz)=0$ for
      $1\le i\le d$} \}.
   \end{displaymath}

   \textnormal{(c)} If $A=\kk\aal\val^{-1}(\D)\aar$, then
   $\Max(A)$ is identified with the set of
   $G$-orbits on
   \begin{displaymath}
      \val^{-1}(\D) \stackrel{\textnormal{def}}{=}\{\bz\in(\kka^\times)^d:
      \val(\bz)\in\D\}.
   \end{displaymath}

   \textnormal{(d)} If
   $A=\kk\aal\val^{-1}(\D)\aar/\<f_1,\dots,f_r\>$, then $\Max(A)$
   is identified with the set of $G$-orbits on
   \begin{displaymath}
      \{\bz\in\val^{-1}(\D):f_i(\bz)=0 \text{ for $1\le i\le d$}\}.
   \end{displaymath}
\end{proposition}

\begin{proof}
   (a) is standard, (b) follows from Proposition
   \ref{prop:tate-algebra}, (c) follows from (b) and the proof of
   Proposition \ref{prop:affinoid-proof}, and (d) follows from (c).
\end{proof}

\subsection{Proof of Theorem
   \ref{thm:main}}\label{subsec:main-proof}

Let $X\subset\GGm^d$ be an irreducible variety. We have already
observed that $\ST(X)\subset\BG(X)$. Since $\BG(X)$ is a
$\Gamma$-rational polyhedral set, we may reduce, as in the proof
of Theorem \ref{thm:A(f)-tropical}, to the following case.

\begin{lemma}\label{lem:BG-proof}
   If $\mathbf{0}\in\BG(X)$, then $\mathbf{0}\in\ST(X)$.
\end{lemma}

\begin{proof}
   Suppose that $\mathbf{0}\notin\ST(X)$. Then 
   $X(\kka)\cap\val^{-1}(\mathbf{0})=\emptyset$. Let $X$ be
   defined by the polynomials $f_1$, $\dots$, $f_r$. Put
   $\widehat{A}=\kk\aal
   \val^{-1}(\bzero)\aar/\<f_1,\dots,f_r\>$. By Proposition~
   \ref{prop:maximal}(d),
   \begin{displaymath}
      \Max(\widehat{A})=\bigl( X(\kka)\cap\val^{-1}(\bzero)\bigr)
      /G=\emptyset,
   \end{displaymath}
   so that $\widehat{A}=0$. Hence there are $g_1$, $\dots$,
   $g_r\in \kk\aal \val^{-1}(\bzero)\aar$ such that
   $1=f_1g_1+\dots+f_rg_r$.
   
   On the other hand, suppose that $\bzero\in\BG(X)$. Then there
   is a valuation $w$ on
   $A=\kk[x_1^{\pm1},\dots,x_d^{\pm1}]/\<f_1,\dots,f_r\>=\kk[X]$
   such that $w(x_i)=0$ for $1\le i \le d$. We can consider $w$ as
   a ring valuation on $\kk[x_1^{\pm1},\dots,x_d^{\pm1}]$ that
   equals $\infty$ on $\<f_1,\dots,f_r\>$. Since $w(x_i)=0$, we
   can extend $w$ by continuity to a valuation $\hatw$ on
   $\kk\aal \val^{-1}(\bzero)\aar$ via
   \begin{displaymath}
      \hatw\Bigl( \sum_{\bn\in\zd} a_{\bn}\bx^{\bn}\Bigr)=
      \lim_{N\to\infty} w\Bigl(\sum_{\bn\in[-N,N]^d}
      a_{\bn}\bx^{\bn} \Bigr).
   \end{displaymath}
   But since $\hatw(f_j)=w(f_j)=\infty$, we obtain that
   \begin{align*}
      0=\hatw(1) &= \hatw(f_1g_1+\dots f_rg_r)
      \ge \min_{j}\{\hatw(f_jg_j)\}\\
      &=\min_{j} \{w(f_j)+\hatw(g_j)\}=\infty.
   \end{align*}
   This contradiction proves the lemma, and hence the theorem.
\end{proof}

\begin{remark}\label{rem:simplify}
   The non-archimedean analytic point of view allows one to
   simplify several proofs in the original Bieri-Groves treatment
   of their sets. Thus, the property of total concavity of
   $\BG(X)$ (see \cite{BG}) follows at once, if formulated for
   $\ST(X)$ instead, from the maximum modulus principle for
   affinoid sets (see Proposition 4 of \cite[Sec.\ 6.2.1]{BGR}).
\end{remark}

\subsection{Reminder on rigid analytic spaces}\label{subsec:reminder} 
   
The basic reference for this section is~ \cite{BGR}, which
contains a complete and accessible treatment of the ideas we use
there. 

Let $A$ be an affinoid
algebra. The set $\Max(A)$ has the following structures (see
\cite[Chap.\ 9]{BGR}):
\begin{enumerate}
  \item[(1)] A Grothendieck topology $\mathscr{T}$ (the strong
    $G$-topology of \cite[Sec.\ 9.1.4]{BGR}), i.e.,
    
    (1a) A family of subsets $U\subset\Max(A)$ called
    \textit{admissible open}, and

    (1b) For any admissible open $U$, a family $\text{Cov}(U)$ of
    coverings of $U$ by admissible open sets contained in $U$,
    satisfying the axioms of \cite[Sec.\ 9.1.1]{BGR}. Such coverings
    are called \textit{admissible}.

   \item[(2)] A sheaf of local rings $\mathscr{O}$ on the above topology
    such that $\mathscr{O}(\Max(A))=A$.
\end{enumerate}

The construction of these objects is given a detailed treatment
in \cite{BGR}. We give here some instructive examples.

\begin{examples}
   (a) If $f_1,\dots,f_r\in A$, then
   \begin{displaymath}
      U_{f_1,\dots,f_r}=\{x\in\Max(A):|f_i(x)|\le1, i=1,\dots,r\}
   \end{displaymath}
   is an admissible open subset called a \textit{Weierstrass domain}
   \cite[Sec.\ 7.2.3]{BGR}.

   (b) If $f_1,\dots,f_r\in A^*$ then the covering of $\Max(A)$ by
   domains of the form
   \begin{displaymath}
      U_{\epsilon_1,\dots,\epsilon_r}=\{x:|f_i(x)|^{\epsilon_i}\le1,
      i=1,\dots,r\}, \quad \epsilon_i=\pm1
   \end{displaymath}
   is admissible. It is called a \textit{Laurent covering}
   \cite[Sec.\ 8.2.2]{BGR}.
\end{examples}

Recall that $G=\Gal(\kka/\kk)$.

\begin{proposition}\label{prop:admissible}
   Let $A=\kk\aal \val^{-1}(\D)\aar$ for a bounded $\G$-rational
   polyhedron $\D$ as in Example \ref{exam:laurent-affinoid}, so that
   $\Max(A)=\val^{-1}(\D)/G$. For any $\G$-rational
   subpolyhedron $\Sigma\subset\D$ the subset
   $\val^{-1}(\Sigma)/G$ is admissible open.
\end{proposition}

\begin{proof}
   Any additional $\G$-rational inequalities defining $\Sigma$
   can be written in the form $|w\cdot z_1^{b_1}\dots
   z_d^{b_d}|\le1$, so the subset in question is a Weierstrass
   domain.
\end{proof}

By definition a (rigid) analytic space over $\kk$ is a system
$Z=(\uz,\ST,\SO_Z)$ consisting of a set $\uz$, a Grothendieck
topology $\ST$ on $\uz$, and a sheaf $\SO_Z$ of rings on $\ST$
such that locally on $\ST$ it is isomorphic to $\Max(A)$ where
$A$ is an affinoid algebra with its Grothendieck topology and
sheaf $\SO$.

\begin{example}\label{exam:analytic-space}
   Every scheme $X$ of finite type over $\kk$ gives rise to an
   analytic space $X^{\an}$ with
   $\underline{X}^{\an}=X(\kka)/G$ (see
   \cite[Sec.\ 9.3.4]{BGR}). We are particularly interested in
   the case when $X\subset\GGm^d$ is a closed subscheme. In this
   case for each bounded $\G$-rational polyhedron $\D\subset\rd$
   the intersection
   \begin{displaymath}
      \bigl(X(\kka)\cap\val^{-1}(\D)\bigr)/G
   \end{displaymath}
   is an admissible subset in $\ux^{\an}$.
\end{example}

\subsection{Connectedness and irreducibility for analytic spaces.
Proof of Theorem \ref{thm:connected}} A rigid analytic space $Z$
is called \textit{disconnected} if it admits an admissible open
covering consisting of two disjoint subspaces. Otherwise $Z$ is
called \textit{connected}. One says that $Z$ is
\textit{irreducible} if the normalization of $Z$ (see \cite[Sec.\ 
2.1]{Co}) is connected. In particular an irreducible analytic
space is connected. The following result of Conrad will be
crucial for us.

\begin{theorem}[Thm\ 2.3.1 of \cite{Co}]\label{thm:conrad}
   Let $X$ be an irreducible algebraic variety over $\kk$. Then
   the analytic space $X^{\an}$ is irreducible and, in
   particular, is connected.
\end{theorem}

We now assume that $X$ is a closed subvariety in $\GGm^d$, and let
$\ST=\ST(X)\subset \rd$ be its tropical variety. Suppose that
$\ST$ is disconnected: $\ST=\SB\sqcup\SC$ where $\SB$ and $\SC$
are open and closed in $\ST$. Define subsets
$B,C\subset\ux^{\an}=X(\kka)/G$ by
\begin{displaymath}
   B=(X(\kka)\cap\val^{-1}(\SB))/G
\end{displaymath}
and similarly for $C$. Then, clearly, $\ux^{\an}=B\sqcup C$. To
establish Theorem \ref{thm:connected} it is therefore enough to
prove the following.

\begin{proposition}\label{prop:admissible}
   \textnormal{(a)} $B$ and $C$ are admissible open in
   $\ux^{\an}$.

   \textnormal{(b)} The covering of $\ux^{\an}$ by $B$ and $C$ is
   admissible.
\end{proposition}

\begin{proof}
   We recall two properties which hold for the Grothendieck
   topology of any rigid analytic space $Z$ (\cite[p.\
   339]{BGR}).

   (G1) Let $U$ be admissible open in $\uz$ and $V\subset U$ be a
   subset. Assume there exists an admissible covering $\{U_i\}$
   of $U$ such that $V\cap U_i$ is admissible open in $\uz$ for
   all $i$. Then $V$ is admissible open in $\uz$.

   (G2) Let $\mathscr{U}=\{U_i\}_{i\in I}$ be a covering of an
   admissible open $U\subset \uz$ such that $U_i$ is admissible
   open in $\uz$ for each $i$. Assume that $\mathscr{U}$ has a
   refinement which is an admissible covering of $U$. Then
   $\mathscr{U}$ itself is an admissible covering of $U$.

   Let $\rd=\bigcup_{i\in I} \D_i$ be a decomposition of $\rd$
   into parallel cubes of sufficiently small size
   $\epsilon$. Then
   $\{\val^{-1}(\D_i)/G:i\in I\}$ is an admissible covering
   of $(\kka^\times)^d/G$ and therefore the sets
   \begin{displaymath}
      D_i=(X(\kka)\cap\val^{-1}(\D_i))/G
   \end{displaymath}
   form an admissible covering of $\ux^{\an}$. Let $J_B=\{i\in
   I:\D_i\cap\SB\ne\emptyset\}$, $J_C=\{i\in
   I:\D_i\cap\SC\ne\emptyset\}$, and $J=J_B\cup J_C$. Since $\SB$
   and $\SC$ are polyhedral, by taking $\epsilon$ small enough,
   we can assume $J_B\cap J_C=\emptyset$. As $D_i=\emptyset$ for
   $i\notin J$, we have $D_i$, $i\in J$ form an admissible
   covering of $\ux^{\an}$. Now, applying (G1) to $Z=U=X^{\an}$,
   $V=B$,
   $U_i=D_i$, $i\in J$ we find that $V\cap U_i=U_i$ for $i\in
   J_B$, and $V\cap U_i=\emptyset$ for $i\notin J_B$, so $V=B$ is
   admissible open. Similarly for $C$. This proves part (a) of
   the proposition. Part (b) follows at once from (G2), as the
   covering $\ux^{\an}=B\cup C$ has a refinement
   $\ux=\bigcup_{i\in J} D_i$ which is admissible.
\end{proof}

\section{Adelic amoebas in algebra and
dynamics}\label{sec:connections}\marginpar

In this section we briefly describe two situations in
which adelic amoebas have already implicitly appeared.

\subsection{The Bieri-Strebel geometric invariant}

Let $R$ be a commutative ring with~ $1$, and $\Rx$ denote the ring
$R[x_1^{\pm1},\dots,x_d^{\pm1}]$ of Laurent polynomials over
$R$. Let $\Sd$ denote the unit sphere in $\rd$. For each
$\bu\in\Sd$ let $H_{\bu}=\{\bv\in\rd:\bu\cdot\bv\le0\}$ be the
half-space with outward normal $\bu$. There is a continuum of
subrings $R_{\bu}=R[\bx^{\bn}:\bn\in\zd\cap H_{\bu}]$ of $R$ as
$\bu$ varies over $\Sd$.

Suppose that $M$ is a finitely generated $\Rx$-module. For which
$\bu\in\Sd$ does $M$ remain finitely generated over
$R_{\bu}$? This question led Bieri and Strebel \cite{BS} to
define their geometric invariant for
 $R[\xpm]$-modules $M$ as
\begin{displaymath}
   \Sigma_M^c =\{\, \bu\in\Sd:M \text{ is \textit{not} finitely
   generated over $R_{-\bu}$} \,\}
\end{displaymath}
(we have used a negative sign in $R_{-\bu}$ since Bieri-Strebel
use inward normals).  This invariant has proved crucial in
answering a number of important algebraic questions. For example,
they show that certain $R[\xpm]$-modules $M$ are finitely
presented if and only if $\Sigma_M^c$ does not contain any pair
of antipodal points (a condition reminiscent of the dynamical
notion of totally non-symplectic).

The geometric invariant $\Sigma_M^c$ can be obtained from
Bieri-Groves sets as follows. Let $v$ be a (ring) valuation on
$R$. For an ideal $\fa$ in $\Rx$ define $\SW(\Rx/\fa)$ to be the
set of all valuations $w$ on $\Rx/\fa$ extending $v$ such that
$w(x_i)<\infty$ for $1\le i\le d$. Define the Bieri-Groves set to
be 
\begin{displaymath}
   \BG_v(\Rx/\fa)=\bigl\{\, \bigl(w(x_1),\dots,w(x_d)\bigr):
   w\in \SW(\Rx/\fa)\, \bigr\}.
\end{displaymath}
Let $\rho_0\colon\rd\smallsetminus\{\bzero\}\to\Sd$ be radial
projection, and extend $\rho_0$ to subsets of $\rd$ by 
$\rho(E)=\rho_0(E\smallsetminus\{\bzero\})$. Then Bieri and
Groves \cite[Thm.\ 8.1]{BG} showed the following.

\begin{theorem}\label{thm:sigma-set}
   Let $M$ be a Noetherian $\Rx$-module and $\fa$ be the
   annihilator of $M$. Then
   \begin{displaymath}
      \Sigma_M^c=\bigcup_{v(R)\ge0}
      \rho\bigl(\BG_v(\Rx/\fa)\bigr),
   \end{displaymath}
   where the union is over all valuations on $R$ which are
   nonnegative.
\end{theorem}

In the special case $R=\ZZ$ the geometric
invariant $\Sigma_M^c$ is related to the part of the adelic amoeba
that corresponds to the finite primes $p$.

\begin{corollary}\label{cor:sigma-set}
   Let $M$ be a Noetherian $\Zx$-module and $\fa$ be the
   annihilator of $M$. Suppose $M$ (or equivalently $\Rx/\fa$)
   is torsion-free
   as a module over $\ZZ$, and define $X\subset\GGm^d$ to be the
   algebraic variety defined by $\fa$ over the field $\QQ$.
   Then
   \begin{displaymath}
      \Sigma_M^c=\bigcup_{p<\infty}
      \rho\bigl(\ST(X(\Qpbar))\bigr),
   \end{displaymath}
   where the union is over all rational prime numbers,
   and $X(\Qpbar)$ is the variety defined by $\fa$
   over the algebraic closure $\Qpbar$ of the $p$-adic rationals.
\end{corollary}

In other words the radial projections of the negatives of the
$p$-adic amoebas of $\fa$ describe $\Sigma_M^c$ completely. Since
there are only finitely many distinct $p$-adic amoebas, this is
actually a finite union. An explicit algorithm for computing this
union, using universal Gr\"obner bases and Fitting ideals, is
described in Proposition~ 6.6 of~ \cite{ELMW}.

\begin{proof}
   By Theorem \ref{thm:sigma-set} $\Sigma_M^c$ can be calculated
   via the Bieri-Groves sets. By Theorem \ref{thm:main} we know
   that $\BG(X)=\ST(X)$ over every fixed field $\kk$ with some
   fixed valuation $v$. However, Theorem \ref{thm:sigma-set} uses
   ring valuations $w$, that are allowed to have $w(n)=\infty$
   for nonzero $n\in\ZZ$.  For the corollary we need to show that
   the restriction to the $p$-adic valuations does not change the
   statement (under the torsion-free assumption).
   
   So assume that $w$ is a valuation on $\Zx/\fa$ with
   $\rho(w(x_1),\ldots,w(x_d))=\bv\ne\bzero$, and $w(p)=\infty$
   for some prime number $p$. We claim that the direction $\bv$
   is also captured by the $p$-adic valuation.
   
   Suppose $\fp_1,\ldots,\fp_k$ are the associated prime ideals
   to $\Zx/\fa$. Then $w(\fp_i)=\{\infty\}$ for some $i$.
   Let $\fp=\fp_i$, and let $\FF$ be the field of fractions of
   $\Zx/\fp$. Note that $\FF$ has characteristic zero by
   assumption. Then \cite[Thm.~C2]{BG} describes
   $\BG_p(\Zx/\fp)$ (using the $p$-adic valuation) near the
   origin, and implies that $\{r\bv:r\in
   [0,\epsilon]\}\subset\BG_p(\Zx/\fp)$ for some
   $\epsilon>0$.  Since
   $\BG_p(\Zx/\fp)\subset\ST(X(\Qpbar))$, the claim
   follows.
\end{proof}

\subsection{Expansive subdynamics of algebraic actions} 
We now turn to dynamics. Again let $M$ be a module over $\Zx$.
There is a corresponding \textit{algebraic $\zd$-action} $\alm$
on a compact group $X_M$ defined as follows. Consider $M$ as a
discrete abelian group, and let $X_M$ be its compact Pontryagin
dual group. For $\bn\in\zd$ let $\alm^{\bn}$ be the automorphism
of $X_M$ dual to the automorphism of $M$ given by multiplication
by $\bx^{\bn}$. This process can be reversed, so that given an
algebraic $\zd$-action by automorphisms of a compact abelian
group, there is a corresponding $\Zx$-module via duality. The
book \cite{Sch} contains all necessary background and a wealth of
examples.

A framework for studying general topological $\zd$-actions was
developed in \cite{BL}, focusing on the key idea of expansiveness
along half-spaces. Fix a metric $\delta$ on $X_M$ compatible with
its topology. Then $\alm$ is called \textit{expansive along
$H_\bu$} if there is an $\epsilon>0$ so that if $\xi$ and $\eta$
are two points in $X_M$ with
$\delta\bigl(\alm^{\bn}(\xi),\alm^{\bn}(\eta)\bigr)<\epsilon$ for
all $\bn\in H_{\bu}\cap\ZZ^d$, then $\xi=\eta$. The
\textit{nonexpansive set of $\alm$} is
\begin{displaymath}
   \N(\alm)=\{\,\bu\in\Sphd:\text{ $\alm$ is \textit{not}
   expansive along $H_{\bu}$}\,\}.
\end{displaymath}

This set turns out to be closed in $\Sphd$. The expansive
subdynamics philosophy advocated in \cite{BL} says that dynamical
properties of $\alm$ restricted to subspaces are either constant
or vary nicely within a connected component of the complement of
$\N(\alm)$, but typically change abruptly when passing from one
connected component to another, analogous to a phase transition.
The description of lower dimensional entropy in \cite[Sec.\ 6]{BL}
is an example of this philosophy in action.

It is therefore natural to ask for an explicit calculation for
the nonexpansive set for an algebraic
$\zd$-action. This was done in \cite[Prop.~4.9]{ELMW} using the
complex amoeba and $\Sigma_M^c$. Combining this with Corollary
\ref{cor:sigma-set} shows that the nonexpansive set is the radial
projection of an adelic amoeba.

\begin{theorem}\label{thm:expansive}
   Let $M$ be a Noetherian $\Zx$-module and $\fa$ be its
   annihilator. Suppose that $X_M$ is connected (or, equivalently,
   that $\Zx/\fa$ is torsion-free over~ $\ZZ$). Then
   \begin{displaymath}
      \N(\alm)=\bigcup_{p\le\infty}
      \rho\bigl(\SA_{\Qpbar}(\fa)\bigr)
      =\rho\bigl(\SA_{\AA}(\fa)\bigr).
   \end{displaymath}
\end{theorem}

\begin{figure}[htbp]
   \begin{center}
      \scalebox{1}{\includegraphics{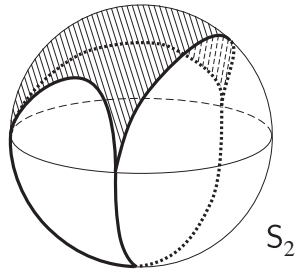}}
      \caption{The nonexpansive set of $\<1+x+y,z-2\>$}
      \label{fig:space-helmet}
   \end{center}
\end{figure}

\begin{example}\label{exam:space-helmet}
   Let $d=3$, $\fa=\<1+x+y,z-2\>$, and $M=\Zx/\fa$. Then
   $\N(\alm)$ is depicted in Figure \ref{fig:space-helmet}. The
   portion above the equator is the radial projection of the
   complex amoeba of $\fa$, the part below the equator is the
   radial projection of the $2$-adic amoeba, and the three points
   on the equator come from the $p$-adic amoeba for $p\ne2$
   (which also form Bergman's logarithmic limit set
   \cite{Bergman} of $\Var(\fa)$). Here the expansive components
   are the three lobes of $\mathsf{S}_2$ in the complement of
   $\N(\alm)$. It is perhaps interesting to note that entropy
   considerations show that none of these components can contain
   a pair of antipodal points.
\end{example}

\end{document}